\documentclass[12pt,eqno,]{article}
\usepackage{amsmath,amsthm,amssymb}

\makeatletter

\newcommand{\Rmnum}[1]{\expandafter\@slowromancap\romannumeral #1@}
\makeatother

\title{\large Optimal two parameter bounds for the Seiffert mean}
\author{\small  Yu-Ming Chu, Miao-Kun Wang and Ye-Fang Qiu}
\date{}

\begin{document}
\maketitle

\begin{center}
{\footnotesize\rm
\noindent Department of Mathematics, Huzhou Teachers College, Huzhou 313000, China\\
Correspondence should be addressed to Yu-Ming Chu,
chuyuming@hutc.zj.cn}
\end{center}

\bigskip \noindent{\bf Abstract}: In this note we obtain sharp
bounds for the Seiffert mean in terms of a two parameter family of
means. Our results generalize and extend the recent bounds presented
in the Journal of Inequalities and Applications (2012) and Abstract
and Applied Analysis (2012).

\noindent{\bf 2010 Mathematics Subject Classification}: 26E60.

\noindent{\bf Keywords}: Seiffert mean, root mean square,
contraharmonic mean.

\section{Introduction}
\hspace{4mm} For $a,b>0$ with $a\neq b$, the Seiffert mean $T(a,b)$,
root mean square $S(a,b)$ and contraharmonic mean $C(a,b)$ are
defined by
\begin{equation}
T(a,b)=\frac{a-b}{2\arctan[(a-b)/(a+b)]},
\end{equation}
\begin{equation}
S(a,b)=\sqrt{\frac{a^2+b^2}{2}}
\end{equation}
and
\begin{equation}
C(a,b)=\frac{a^2+b^2}{a+b},
\end{equation}
respectively. It is well known that the inequalities
\begin{equation*}
T(a,b)<S(a,b)<C(a,b)
\end{equation*}
hold for all  $a,b>0$ with $a\neq b$ .

Recently, $T(a,b)$, $S(a,b)$ and $C(a,b)$ have been the subject of
intensive research. In particular, many remarkable inequalities and
properties for these means can be found in the literature [1-8].

For $\alpha,\beta,\lambda,\mu\in(1/2,1)$, very recently Chu et al.
[9, 10] proved that the inequalities
\begin{equation}
S(\alpha a+(1-\alpha)b,\alpha b+(1-\alpha)a)<T(a,b)<S(\beta
a+(1-\beta)b,\beta b+(1-\beta)a)
\end{equation}
and
\begin{equation}
C(\lambda a+(1-\lambda)b,\lambda b+(1-\lambda)a)<T(a,b)<C(\mu
a+(1-\mu)b,\mu b+(1-\mu)a)
\end{equation}
hold for all $a,b>0$ with $a\neq b$ if and only if $\alpha \leq
(1+\sqrt{16/\pi^2-1})/2$, $\beta\geq(3+\sqrt{6})/6$, $\lambda\leq
(1+\sqrt{4/\pi-1})/2$ and $\mu\geq(3+\sqrt{3})/6$.

Let $t\in(1/2,1)$, $p\geq 1/2$ and
\begin{equation}
Q_{t,p}(a,b)=C^{p}(ta+(1-t)b, tb+(1-t)a)A^{1-p}(a,b),
\end{equation}
where $A(a,b)=(a+b)/2$ is the classical arithmetic mean of $a$ and
$b$. Then from (1.2), (1.3) and (1.6) we clearly see that
\begin{equation*}
Q_{t,1/2}(a,b)=S(ta+(1-t)b, tb+(1-t)a),
\end{equation*}
\begin{equation*}
Q_{t,1}(a,b)=C(ta+(1-t)b, tb+(1-t)a)
\end{equation*}
and $Q_{t,p}(a,b)$ is strictly increasing with respect to
$t\in(1/2,1)$ for fixed $a,b>0$ with $a\neq b$.

It is natural to ask what are the greatest value $t_{1}=t_{1}(p)$
and the least value $t_{2}=t_{2}(p)$ in $(1/2,1)$ such that the
double inequality
\begin{equation*}
Q_{t_{1},p}(a,b)<T(a,b)<Q_{t_{2},p}(a,b)
\end{equation*}
holds for all $a,b>0$ with $a\neq b$ and $p\geq 1/2$. The aim of
this paper is to answer this question, our main result is the
following Theorem 1.1.

\medskip
{\bf Theorem 1.1.} If $t_{1},t_{2}\in(1/2,1)$ and
$p\in[1/2,\infty)$, then the double inequality
\begin{equation}
Q_{t_{1},p}(a,b)<T(a,b)<Q_{t_{2},p}(a,b)
\end{equation}
holds for all $a,b>0$ with $a\neq b$ if and only if $t_{1}\leq
1/2+[\sqrt{(4/\pi)^{1/p}-1}]/2$ and $t_{2}\geq 1/2+\sqrt{3p}/(6p)$.

\medskip
{\bf Remark 1.1.} If we take $p=1/2$ and $p=1$ in Theorem 1.1, then
inequality (1.7) reduces to inequalities (1.4) and (1.5),
respectively.

\section{Proof of Theorem 1.1}
\hspace{4mm} \setcounter{equation}{0} In order to prove Theorem 1.1
we need two lemmas, which we present in this section.

\medskip
{\bf Lemma 2.1.} (see [11, Theorem 1.25]). For $-\infty<a<b<\infty$,
let $f,g:[a,b]\rightarrow{\mathbb{R}}$ be continuous on $[a,b]$, and
be differentiable on $(a,b)$, let $g'(x)\neq 0$ on $(a,b)$. If
$f^{\prime}(x)/g^{\prime}(x)$ is increasing (decreasing) on $(a,b)$,
then so are
$$\frac{f(x)-f(a)}{g(x)-g(a)}\ \ \mbox{and}\ \ \frac{f(x)-f(b)}{g(x)-g(b)}.$$
If $f^{\prime}(x)/g^{\prime}(x)$ is strictly monotone, then the
monotonicity in the conclusion is also strict.

\medskip {\bf Lemma 2.2.} Let $u\in[0,1]$, $p\geq 1/2$ and
\begin{equation}
f_{u,p}(x)=p\log(1+ux^2)-\log{x}+\log{\arctan{x}}.
\end{equation}
Then

(1) $f_{u,p}(x)>0$ for $x\in(0,1)$ if and only if $3pu\geq 1$;

(2) $f_{u,p}(x)<0$ for $x\in(0,1)$ if and only if $1+u\leq
(4/\pi)^{1/p}$.

\medskip
{\bf Proof.} By (2.1) and simple computations one has
\begin{equation}
\lim_{x\rightarrow 0}f_{u,p}(x)=0,
\end{equation}
\begin{equation}
\lim_{x\rightarrow
1}f_{u,p}(x)=p\log(1+u)+\log\left(\frac{\pi}{4}\right)
\end{equation}
and
\begin{align}
{f_{u,p}}'(x)=&\frac{2pu x}{1+u x^2}+\frac{1}{(1+x^2)\arctan{x}}-\frac{1}{x}\nonumber\\
=&\frac{u[(2p-1)x^2(1+x^2)\arctan{x}+x^3]-[(1+x^2)\arctan{x}-x]}{x(1+x^2)(1+ux^2)\arctan{x}}\nonumber\\
=&\frac{(2p-1)x^2(1+x^2)\arctan{x}+x^3}{x(1+x^2)(1+ux^2)\arctan{x}}[u-g(x)],
\end{align}
where
\begin{equation*}
g(x)=\frac{(1+x^2)\arctan{x}-x}{(2p-1)x^2(1+x^2)\arctan{x}+x^3}.
\end{equation*}

Let $g_{1}(x)=\arctan{x}-x/(1+x^2)$, and
$g_{2}(x)=(2p-1)x^2\arctan{x}+x^3/(1+x^2)$. Then
\begin{equation}
g(x)=\frac{g_{1}(x)}{g_{2}(x)}, \quad g_{1}(0)=g_{2}(0)=0
\end{equation}
and
\begin{equation}
\frac{{g_{1}}'(x)}{{g_{2}}'(x)}=\frac{1}{(2p-1)[(1+x^2)^2\arctan{x}]/x+px^2+p+1}.
\end{equation}

It is not difficult to verify that the function $x\rightarrow
[(1+x^2)^2\arctan{x}]/x$ is strictly increasing from $(0,1)$ onto
$(1,\pi)$, hence (2.6) implies that ${g_{1}}'(x)/{g_{2}}'(x)$ is
strictly decreasing in $(0,1)$. Therefore, $g(x)$ is strictly
decreasing in $(0,1)$ follows from Lemma 2.1 and (2.5) together with
the monotonicity of ${g_{1}}'(x)/{g_{2}}'(x)$. Moreover, making use
of  l'H\^{o}pital's rule we get
\begin{equation}
\lim\limits_{x\rightarrow 0}g(x)=\frac{1}{3p}
\end{equation}
and
\begin{equation}
\lim\limits_{x\rightarrow 1}g(x)=\frac{\pi-2}{(2p-1)\pi+2}.
\end{equation}

We divide the proof into three cases.

{\bf Case 1} $u\geq 1/(3p)$. Then from (2.4) and (2.7) together with
the monotonicity of $g(x)$ lead to conclusion that $f_{u,p}(x)$ is
strictly increasing in $(0,1)$. Therefore $f_{u,p}(x)>0$ for
$x\in(0,1)$ follows from (2.2) and the monotonicity of $f_{u,p}(x)$.

{\bf Case 2} $u\leq (\pi-2)/[(2p-1)\pi+2]$. Then from (2.4) and
(2.8) together with the monotonicity of $g(x)$ we clearly see that
$f_{u,p}(x)$ is strictly decreasing in $(0,1)$. Therefore
$f_{u,p}(x)<0$ for $x\in(0,1)$ follows from (2.2) and the
monotonicity of $f_{u,p}(x)$.

{\bf Case 3} $(\pi-2)/[(2p-1)\pi+2]<u<1/(3p)$. Then from (2.4),
(2.7) and (2.8) together with the monotonicity of $g(x)$ we know
that there exists $x_{0}\in(0,1)$ such that $f_{u,p}(x)$ is strictly
decreasing in $(0,x_{0})$ and strictly increasing in $(x_{0},1)$.

Let $h_{p}(u)=\lim\limits_{x\rightarrow 1}f_{u,p}(x)$. Then it
follows from (2.3) that
\begin{equation}
h_{p}(u)=p\log(1+u)+\log\left(\frac{\pi}{4}\right).
\end{equation}

Cases 1 and 2 implies that
\begin{equation}
h_{p}(\frac{1}{3p})=p\log\left(1+\frac{1}{3p}\right)+\log\left(\frac{\pi}{4}\right)>0
\end{equation}
and
\begin{equation}
h_{p}\left(\frac{\pi-2}{(2p-1)\pi+2}\right)=p\log\left[1+\frac{\pi-2}{(2p-1)\pi+2}\right]+\log\left(\frac{\pi}{4}\right)<0.
\end{equation}

From (2.9) we clearly see that $h_{p}(u)$ is strictly increasing in
$[(\pi-2)/[(2p-1)\pi+2],1/(3p)]$, then (2.10) and (2.11) lead to
conclusion that there exists
$u_{0}=(4/\pi)^{1/p}-1\in((\pi-2)/[(2p-1)\pi+2],1/(3p))$ such that
$h_{p}(u)<0$ for $u\in [(\pi-2)/[(2p-1)\pi+2],u_{0})$ and
$h_{p}(u)>0$ for $u\in (u_{0},1/(3p)]$, where
$u_{0}=(4/\pi)^{1/p}-1$ is the unique solution of the equation
$h_{p}(u)=0$ in $[(\pi-2)/[(2p-1)\pi+2],1/(3p)]$. Therefore,
$f_{u,p}(x)<0$ for all $x\in(0,1)$ if and only if
$(\pi-2)/[(2p-1)\pi+2]<u<u_{0}=(4/\pi)^{1/p}-1$ follows from (2.2),
(2.3) and (2.9) together with the piecewise monotonicity of
$f_{u,p}(x)$.

\medskip
{\bf Proof of Theorem 1.1.} Since both $Q_{t,p}(a,b)$ and $T(a,b)$
are symmetric and homogeneous of degree 1. Without loss of
generality, we assume that $a>b$. Let $x=(a-b)/(a+b)\in(0,1)$. Then
from (1.1) and (1.6) we get
\begin{align}
&\log\left(\frac{Q_{t,p}(a,b)}{T(a,b)}\right)=\log\left(\frac{Q_{t,p}(a,b)}{A(a,b)}\right)-\log\left(\frac{T(a,b)}{A(a,b)}\right)\nonumber\\
&=p\log\left[1+(1-2t)^2x^2\right]-\log{x}+\log{\arctan{x}}.
\end{align}

Therefore, Theorem 1.1 follows from Lemma 2.2 and (2.12).\\

{\bf Acknowledgement:} This work was supported by the Natural
Science Foundation of China (Grant Nos. 11071059, 11071069,
11171307), and the Innovation Team Foundation of the Department of
Education of Zhejiang Province (Grant no. T200924).

\end{document}